\documentclass{amsart}[12pt]
\usepackage{amsmath,amssymb,amscd}

\DeclareMathOperator{\reg}{reg}

\DeclareMathOperator{\sat}{sat}
\DeclareMathOperator{\Ht}{ht}
\newtheorem{thm}{Theorem}[section]

\newtheorem{lem}[thm]{Lemma}
\newtheorem{prop}[thm]{Proposition}
\newtheorem{defin}[thm]{Definition}

\newcommand{\p}[1]{\ensuremath{\mathbb{P}^#1}}

\newcommand{\coh}[3]{\ensuremath{\mathrm{H}^{#1}(#2, #3) } }

\begin{document}

\title[Castelnuovo-Mumford regularity
of subspace arrangements]{A sharp bound for the Castelnuovo-Mumford regularity of subspace arrangements}
\author{Harm Derksen}
\address{Harm Derksen, Department of Mathematics, University of Michigan, 525 East University Avenue, Ann Arbor, MI 48109-1109}
\thanks{The research of the first author is partially supported by NSF Grant DMS 0102193.}
\email{hderksen@math.lsa.umich.edu}
\author{Jessica Sidman}
\address{Jessica Sidman, Department of Mathematics, University of Michigan, 525 East University Avenue, Ann Arbor, MI 48109-1109}
\email{jsidman@umich.edu}
\subjclass{Primary 13D02; Secondary 52C35}

\maketitle
\section*{Introduction}
Over the past twenty years rapid advances in computational algebraic geometry have generated increasing amounts of interest in quantifying the ``complexity'' of ideals and modules.  For a finitely generated module $M$ over a polynomial ring $S = k[x_0, \ldots, x_n]$ with $k$ an arbitrary field, we say that $M$ is $r$-regular (in the sense of Castelnuovo and Mumford) if the $i$-th syzygy module of $M$ is generated in degrees less than or equal to $r+i$.  The notion of the regularity of $M$ is key in determining the amount of computational resources that working with $M$ requires.   In this paper we will prove that the ideal of an arrangement of $d$ linear spaces in \p{n} is $d$-regular, answering a question posed by B. Sturmfels.

There is a mysterious gap in our understanding of the behavior of the regularity of ideals.  Namely, it is known that in the general case the regularity of an ideal may be as bad as doubly exponential in the degrees of its minimal generators and the number of variables in the ambient ring, and that this is essentially the worst case \cite{mayr-meyer}, \cite{bayer-stillman88}, \cite{giusti}, \cite{yap}, \cite{galligo}.  However, for an ideal $I$ corresponding to a smooth variety in characteristic zero, the regularity of $I$ is linear in terms of geometric data (\cite{gruson-lazarsfeld-peskine}, \cite{pinkham}, \cite{lazarsfeld87}, \cite{ran}, \cite{bertram-ein-lazarsfeld},  \cite{bayer-mumford}).  From this perspective, the ideals of subspace arrangements provide a test case for the behavior of the regularity of ideals of reduced but possibly singular schemes.  

Ideals of subspace arrangements also have connections with invariant theory and combinatorics.  In \cite{noether}, E. Noether gave a proof which showed that for the action of a finite group $G$ on a polynomial ring $k[V]$ induced by a representation of $G$ on $V$, the generators of $k[V]^G$ have degree less than or equal to the order of $G$ as long as the characteristic of $k$ is zero or if the characteristic is positive and greater than the order of $G$.  For some time it was an open question as to whether or not the same result was true in positive characteristic with the stipulation that the characteristic of $k$ not divide the order of $G$.  Although the statement was recently shown to be true  (see \cite{fleischmann}, \cite{fogarty}, \cite{derksen-kemper} for details), our main theorem yields an additional proof completing the circle of ideas introduced by the first author in \cite{derksen}.  Additionally, the ideals of subspace arrangements are studied by combinatorialists in relation to the computation of Betti numbers of the complement of an arrangement of linear subspaces in affine space \cite{montaner-lopez-armengou}.

The ideals of subspace arrangements have an appealing inductive structure and are algebraically closely related to the products of the linear ideals defining the constituent subspaces of the arrangement.  The second author has used sheaf-theoretic techniques to study subspace arrangements of low dimension \cite{sidman}.  The proof given here is algebraic in nature and  illustrates the close relationship of the intersection of linear ideals to their product as it relies heavily on ideas introduced by A. Conca and J. Herzog to control the regularity of the product of linear ideals in \cite{conca-herzog}.  We recall basic notions related to regularity in \S 1 and present the proof of our result in \S 2.

We are very grateful to Bernd Sturmfels for suggesting the collaboration and Rob Lazarsfeld for his support, encouragement, and guidance.  We are also deeply indebted  to Aldo Conca and J{\"u}rgen Herzog for their many kind communications and in particular for sending us early versions of their preprints.  Additionally, we thank Mel Hochster and Tamon Stephen for help searching the literature in commutative algebra and combinatorics, respectively, Dan Rogalski for helpful discussions, and Lawrence Ein and Karen Smith for their comments on the proof.

\section{Basic notions and conventions}

Let $S = k[x_0, \ldots, x_n]$ where $k$ is an infinite field of arbitrary characteristic and ${\mathfrak m} = (x_0, \ldots, x_n)$.  In this section we recall basic definitions and facts concerning regularity.

Of the various formulations of the definition of regularity, the following is the most illuminating from a computational point of view.

\begin{defin}\label{defin: regdef}
Let $M$ be a finitely generated graded module over $S$.  Take a minimal resolution of $M$ by free graded $S$-modules:
\[0 \longrightarrow \oplus^{r_n}_{\alpha = 1}S(-e_{\alpha, n}) \longrightarrow \ldots \longrightarrow \oplus^{r_0}_{\alpha = 1} S(- e_{\alpha, 0}) \longrightarrow M \longrightarrow 0\] We say that $M$ is $r$-regular if $\deg(e_{\alpha, i}) \leq r+i $ for all $\alpha, i$ and that the regularity of $M$, denoted $\reg(M)$, is the least integer for which this holds.
\end{defin}
One can also formulate the definition in terms of local cohomology and sheaf cohomology.  The following theorem (Definition 3.2 in \cite{bayer-mumford}) establishes the relationship between the regularity of an ideal and its associated sheaf. 

\begin{thm}\label{thm: regdef}
Let $I$ be a homogeneous ideal of $S$ and let $\widetilde{I}$ be the associated coherent sheaf of ideals.  Then the following conditions are equivalent.
\begin{itemize}
\item[(a)] The natural map $I_r \longrightarrow \coh{0}{\p{n}}{\widetilde{I}(r)}$ is an isomorphism and \coh{i}{\p{n}}{\widetilde{I}(r-i)} = 0, $1 \leq i \leq n.$

\item[(b)] The natural maps  $I_d \longrightarrow \coh{0}{\p{n}}{\widetilde{I}(d)}$ are isomorphisms for all $d \geq r$ and  \coh{i}{\p{n}}{\widetilde{I}(d)} = 0, $d + i \geq r$, $i \geq 1.$

\item[(c)] The ideal $I$ is $r$-regular.
\end{itemize}
\end{thm}

A concept that arises frequently in the study of regularity is:

\begin{defin}  Let $I$ be a homogeneous ideal of $S$.  Then the saturation of $I$, denoted $I^{\sat},$ is \[I^{\sat} := \{s \in S: s{\mathfrak m}^k \subseteq I, k \gg 0\}.\]  The saturation degree of $I$, denoted $\sat(I)$, is the least integer $k$ for which $I$ agrees with $I^{\sat}$ in degrees $k$ and higher.
\end{defin}

It is also important to understand how regularity behaves in short exact sequences.

\begin{lem} [Corollary 20.19 in \cite{eisenbud}]\label{lem: seq}
Let $A$, $B$, and $C$ be finitely generated graded modules over $S$ and let \[ 0 \longrightarrow A \longrightarrow B \longrightarrow C \longrightarrow 0\] be a short exact sequence.  Then
\begin{itemize}
\item[(a)]  $\reg(A) \leq \max\{\reg(B), \reg(C) +1\}$.
\item[(b)]  $\reg(B) \leq \max\{\reg(A), \reg(C)\}$.
\item[(c)]  $\reg(C) \leq \max\{\reg(A) -1, \reg(C) \}$.
\end{itemize}
\end{lem}

The lemma below (Lemma 1.8 in \cite{bayer-stillman87}) is the basis for a standard method of working inductively with regularity and is central to our proof in \S 2.  We provide a proof here for the convenience of the reader.  (Readers may also wish to consult \cite{bayer-stillman87} for a proof using local cohomology instead of sheaf cohomology.) 
 
\begin{lem}\label{lem: hyperplane}
Let $I \subseteq k[x_0, \ldots, x_n]$ be a homogeneous ideal, and suppose that $x$ is a linear form that is a nonzerodivisor modulo $I^{\sat}$.  Then \[ \reg(I) = \max\{ \reg(I+(x)), \sat(I)\}.\]
\end{lem}

\begin{proof}
Note that the result will follow if we can show that $I$ is $r$-regular if and only if $I$ is saturated in degrees $r$ and higher and $I+(x)$ is $r$-regular.  It is crucial to observe here that $I \cap (x)$ can be rewritten as $(I:x)x$ and that if $I$ is saturated in degrees $r$  and higher then $(I:x)$ agrees with $I$ in degrees $r$ and higher. 
  
For the forwards implication, assume that $I$ is $r$-regular and note that since $\Gamma_*(\widetilde{I}) = I^{\sat}$, part (b) of Theorem \ref{thm: regdef} implies that an $r$-regular ideal has saturation degree less than or equal to $r$.  To see that $I+(x)$ is $r$-regular, consider the short exact sequence \[ \label{eq: seq} 0 \longrightarrow I \cap (x) \longrightarrow I \oplus (x) \longrightarrow I+(x) \longrightarrow 0.\]  Since $(I:x)$ agrees with $I$ in degrees $r$ and higher, $(I:x)$ is $r$-regular which implies that $(I:x)x$ is $(r+1)$-regular.  The $r$-regularity of $I+(x)$ follows from part (c) of Lemma \ref{lem: seq}.

For the reverse implication note that since $I+(x)$ is $r$-regular, its $r$-th graded piece is isomorphic to $H^0\big{(}\p{n}, (I+(x))^{\sim} (r)\big{)}.$  Furthermore, the $r$-th graded piece of $I \oplus (x)$ injects into $H^0\big{(}\p{n}, (I\oplus(x))^{\sim} (r)\big{)}$ and surjects onto the $r$th graded piece of $I+(x)$.  Therefore, \[H^0\big{(}\p{n}, (I\oplus(x))^{\sim} (r)\big{)} \longrightarrow H^0\big{(}\p{n}, (I+(x))^{\sim} (r)\big{)} \] is also a surjection.  Together with the $r$-regularity of $I+(x)$ this implies that we have isomorphisms \[ \coh{i}{\p{n}}{\widetilde{I}(k-1)} \longrightarrow \coh{i}{\p{n}}{\widetilde{I}(k)} \oplus \coh{i}{\p{n}}{\widetilde{(x)}(k)} \] for all $i > 0$ and $k \ge r-i+1$.   But now the vanishing of cohomology for $k \gg 0$ implies that $\coh{i}{\p{n}}{\widetilde{I}(r-i)} = 0$ for all $i >0$.  Since $I$ is saturated in degrees $r$ and higher by assumption, condition (a) of Theorem \ref{thm: regdef} is satisfied and $I$ is $r$-regular.
\end{proof}

\section{The proof}

We wish to prove the following theorem:

\begin{thm}\label{thm: thetheorem}If $X$ is an arrangement of $d$ linear spaces in \p{n} whose ideals are $I_1, \ldots, I_d$ then $I(X) = I_1 \cap \dots \cap I_d$ is $d$-regular.  
\end{thm}

We should remark here that the statement of Theorem \ref{thm: thetheorem} is optimal.  Indeed, one may take $I$ to be the ideal of a general arrangement of $d$ lines that meet a fixed auxiliary line $L$ in $d$ distinct points. ($L$ does not appear in the arrangement.)  The hypersurface corresponding to any $f \in I$ must meet $L$ in at least $d$ points.  Therefore, if $f$ has degree strictly less than $d$, it must contain the line $L$.  From this we can see that $I$ must have at least one generator of degree at least $d$ and that $\reg(I) \geq d$. 

Note that since $I_1 \cap \dots \cap I_d$ is saturated, our result would follow from Lemma \ref{lem: hyperplane} if we could show that $\reg(I_1 \cap \dots \cap I_d +(x)) \leq d$.  In fact, we prove the stronger statement:

\begin{prop}\label{prop: aux}
Let $I_1, \ldots, I_d$ be the ideals of a linear arrangement in \p{n} and let $L$ be any linear ideal.  Then \[ I_1 \cap \dots \cap I_d + L\] is $d$-regular.
\end{prop}

Observe that $I_1 \cap \dots \cap I_d +{\mathfrak m} = {\mathfrak m}$ and is hence $d$-regular (in fact, 1-regular).  This suggests that we ought to proceed by decreasing induction on $\Ht(L)$ (which is also the vector space dimension of the degree one part of $L$).  The following lemma shows that bounding the saturation degree of $I_1 \cap \dots \cap I_d +L$ is actually the crucial computation:

\begin{lem}
If $I_1 \cap \dots \cap I_d + L'$ is $d$-regular for any linear ideal $L'$ with $\Ht(L') > \Ht(L)$, then the $d$-regularity of $I_1 \cap \dots \cap I_d+L$ will follow if $I_1 \cap \dots \cap I_d + L$ is saturated in degrees $d$ and higher.
\end{lem}

\begin{proof}
If $(I_1 \cap \dots \cap I_d + L)^{\sat} = S$ then it is clear that if $I_1 \cap \dots \cap I_d + L$ is saturated in degrees $d$ and higher then it is $d$-regular.  If $(I_1 \cap \dots \cap I_d + L)^{\sat} \neq S$ then choose a linear form $x$ that is a nonzerodivisor modulo $(I_1 \cap \dots \cap I_d + L)^{\sat}$.  Using Lemma \ref{lem: hyperplane} we have \[ \reg(I_1 \cap \dots \cap I_d + L) = \max\{ \reg(I_1 \cap \dots \cap I_d + L+(x)), \sat(I_1 \cap \dots \cap I_d + L) \}.\]  Since $x$ is a nonzerodivisor modulo $I_1 \cap \dots \cap I_d +L$, it cannot be contained in $L$.  Therefore, $L +(x)$ is a linear ideal with height greater that $\Ht(L)$, so $I_1 \cap \dots \cap I_d + (L + (x))$ is $d$-regular by assumption.
\end{proof}

Thus, the heart of the proof of Proposition \ref{prop: aux} is a statement about saturation which we will control with the method used by Conca and Herzog in \cite{conca-herzog}.

\begin{proof}(Proposition \ref{prop: aux})
 Set \[ J_i := I_1 \cap \dots \cap I_{i-1} \cap I_{i+1} \cap \dots \cap I_d +L\] and \[I:=  I_1 \cap \dots \cap I_d + L.\] Since $I_1 + L$ is 1-regular, we may assume inductively that $J_i$ is $(d-1)$-regular.  

Notice that if $\sum I_i + L \neq {\mathfrak m}$ then we can view the problem as living in a polynomial ring in fewer variables, say $S' = k[x_0, \ldots, x_r]$ with $r < n$.  In this case the generators of the ideal $I = I_1 \cap \dots \cap I_d+L$ also live in $S'$.  If we know that the ideal in $S'$ is $d$-regular, then the $d$-regularity of the ideal in the larger ring follows by the flatness of $S'[x_{r+1}, \ldots, x_n]$ over $S'$.  But we may assume $d$-regularity in the smaller ring inductively since the result is trivial for \p{0}.

Finally, we show that $I$ is saturated in degrees $d$ and higher under the assumption that $\sum I_i +L = {\mathfrak m}$.  Suppose there is a homogeneous form $f$ of degree greater than or equal to $d$ such that $f{\mathfrak m}^s \subseteq I$ for some integer $s > 0$.  Choose a linear form $x$ that is not in $L$ and is a nonzerodivisor modulo all of the $J_i^{sat}$.   Since $\Ht(L+(x)) > \Ht(L)$, by induction we see that $f \in I +(x)$ so $f = h + xf'$ with $h \in I$.  Replace $f$ by $xf'$, and note that $xf'{\mathfrak m}^s \subseteq I \subseteq J_i$ implies that $xf' \in J_i^{sat}$ for each $i$.  Now, since $x$ is a nonzerodivisor modulo $J_i^{sat}$, $xf' \in J_i^{sat}$ implies that $f' \in J_i^{sat}$.  Additionally, we know that $J_i$ is saturated in degrees $d-1$ and higher which implies that $f' \in J_i$.  Write $x = \sum x_i + x'$ where $x_i \in I_i$ and $x' \in L$.  Observe that $x_if'+ x'f' \in I_iJ_i \subseteq I$ for each $i$ which implies that $xf' \in I$ and completes our proof.

\end{proof}

\bibliographystyle{amsplain}

\end{document}